\documentclass[11pt]{article}

\title {First Contact Remarks on Umbra Difference Calculus References Streams}
\author{A.K.Kwa\'sniewski\\  
\\ Higher School of Mathematics and Applied Informatics\\
PL-15-021 Bia{\l}ystok, ul. Kamienna 17, POLAND
\\e-mail: kwandr@uwb.edu.pl}

\usepackage{amsmath,amsthm}

\chardef\bslash=`\\ 
\begin{document}
\maketitle\

\begin{abstract}
The reference links to the modern "classical umbral calculus"
(before that properly called Blissard`s symbolic method) and to
Steffensen -actuarialist.....are numerous. The reference links to
the EFOC (Extended Finite Operator Calculus) founded by  Rota with
numerous outstanding Coworkers, Followers and Others ... and links
to Roman-Rota functional formulation of umbra calculi ... .....are
giant numerous. The reference links to the difference $q$-calculus
- umbra-way treated or without even referring to umbra....
.....are plenty numerous. These reference links now result in
counting  in thousands the relevant papers and many books . The
purpose  of the present attempt  is to offer one of the  keys to
enter the world of those thousands of references. This is the
first glimpse - not much structured - if at all. The place of
entrance was chosen selfish being subordinated to  my present
interests and  workshop purposes. You are welcomed to add your own
information.
\end{abstract}

{\bf Motto} {\bf Herman Weyl}      April 1939 {\it The modern
evolution… has on the whole been marked by a trend of
algebraization} ... - …from an invited address before AMS in
conjunction with the cetential celebration of Duke University.

\section{I. First Contact  Umbral  Remark}

"$R$-calculus" and specifically $q$ -calculus - might be
considered as specific cases of umbral calculus  as illustrated by
a]  and  b] below :\\
 a] Example 2.2.from [1] quotation:...."The example of
$\psi$-derivative:
$$
Q(\partial_{\psi})= R(q\hat Q)\partial_0 \equiv \partial_R
$$
i.e.
$$\Big\{[R(q^n)!]^{-1}\Big \}_{n\geq 0}$$
one may find in [2] where an advanced theory of general quantum
coherent states is being developed. The operator $R(q\hat Q)$ is
not recognized in [2] as an example of $\psi$-derivative".

\vspace{2mm}

b] specifically $q$-calculus - might be considered as a particular
case of umbral calculus .It is so recognized since tenths of years
in numerous references to $q$-calculus ({\bf Roman, Rota} and
Others). From algebraic point of view we deal with specific case
of Extended Umbral Calculus (this being particularly evident in
the "Extended Finite Operator Calculus" formulation ; see:
[3,4,5,1] (http://ii.uwb.edu.pl/akk/publ1.htm ).

\vspace{2mm}

The $q$-case of umbra  - has been -with autonomy  - $q$- evolving
{\bf in a specification dependent way} ...... evolving -with
$q$-autonomy - despite the Ward 1936 paper [6]  "ultimately"
extending the Jackson formulation of $q$-difference calculi to
arbitrary "reasonable" difference calculi - in the spirit of
nowadays Extended Finite Operator Calculus  (see: [1]
math.CO/0312397 ; see there for $\psi$-difference calculus and
also for {\bf Markowsky} [7] general difference calculus -  in
Rota-like operator formulation - compare with  Roman -Rota, Roman
functional formulation).

\vspace{2mm}

We share the conviction with more experienced community that this
historically established autonomy - adds and still may add more
inspiration of the analytical and geometrical character  to the
modern formulation of umbra idea (see: [1] and references therein)
-which is so fruitfully formal algebraic and combinatorial in so
many applications (see: [11-21] -via links - for\textbf{thousands
of  references}).

\vspace{2mm}

{\bf References on modern umbra}

The recent modern umbra was born in 19-th century and is known now
under the name of {\bf Blissard} calculi. For more and references
on this and that see:

[8]Andrew P. Guinand the  survey of elementary mnemonic and
manipulative uses of the unbral method;

[9] Brian D. Taylor: on difference equations via the classical
umbral calculus  where you learn  what classical umbral calculus
is to proceed with next:

[10] G.-C. Rota and B. D. Taylor: just on  the Classical Umbral
Calculus

[11] {\bf Ira M. Gessel} here also apart of applications of the
classical umbral calculus you may learn  what classical umbral
calculus is in the Algebra Universalis Realm. and do not miss
many, {\bf many representative  references therein}.

\vspace{2mm}

Quite a lot of time , years ago:

\vspace{2mm}

[12] {\bf Adams C. R.}   in His "Linear q-difference equations"
contained around {\bf 200} up-that-date (1931) references!

\vspace{2mm}

Quite recently  (19December 2001)

[13] Thomas Ernst in his {\bf Licentiate }elaborate included {\bf
969} references. As natural - the bibliography  in this much
useful [13]  is far from being complete - neither is there
$q$-umbral calculi fully represented.

\vspace{2mm}

Compare with the following  links :

\vspace{2mm}

[14] {\bf Umbral Calculus}
$$http://en2.wikipedia.org/wiki/Umbral_calculusUmbral$$

[15] {\bf Blissard symbolic calculus}:
$$http://www.google.pl/search?q=Blissard+symbolic+calculus$$

[16] A. Di Bucchianico  , D. Loeb    A Selected Survey of Umbral
Calculus w ({\bf 506} references)

[17] {\bf D. Loeb , G. C. Rota}  another survey with quite recent
contributions to the general calculus of finite differences

\vspace{2mm}

see also: $http://arxiv.org/list/math.CO/9502$

\vspace{2mm}

and of course

\vspace{2mm}

[18] {\bf S.Roman} in $ http://www.romanpress.com/me.htm $. Here
is the splendid , competent source from The Source.

\vspace{2mm}

Also  other links below are recommended

[19] http://mathworld.wolfram.com/UmbralCalculus.html

[20] {\bf R. Stanton }: competent source from The Source
$$http://www.math.umn.edu/~stanton/publist.html$$

\vspace{2mm}

altogether with highly appreciate link to  NAVIMA - group $q$
-active also in $q$-research fruitfully

[21]  NAVIMA http://webs.uvigo.es/t10/navima/publications.html

\vspace{2mm}

Though no list is  aposteriori complete  - this presented here and
in statu nascendi  listing of references streams may be actualized
with {\bf your e-mail help}.

\vspace{2mm}

\section{II. Second  Contact Remark } on  $\partial _{q,h}$-{\bf calculus
           -an umbral difference calculus of course}.

$\partial _{q,h}$-difference calculus of

{\bf Hahn} [22] may be reduced to $q$-calculus of {\bf
Thomae-Jackson} [23,24,25] due to the following observation . Let

$$
h \in F ,(E_{q,h}\varphi)(x) =\varphi qx+h)
$$
and  let
\begin{equation}\label{eq1}
(\partial _{q,h}\varphi)(x)= \frac{\varphi(x)-
\varphi(qx+h)}{(1-q)x-h}
\end{equation}

\vspace{2mm}

Then (see Hann [22]and also Appendix A in [26] - fifty years after
Hann)

\vspace{2mm}

\begin{equation}\label{eq2}
\partial_{q,h} = E_{1,{\frac{{-h}}{{1-q}}}} \partial_q E_{1,{\frac
{h}{{1-q}}}} .
\end{equation}

\vspace{2mm}

Due to (2) it is easy now to derive corresponding formulas

including Bernoulli-Taylor $\partial _{q,h}$-formula

obtained  in [27] by the Viskov method [28]  which for
$$
q\rightarrow 1 , h\rightarrow 0
$$ recovers the content of one of
the examples in [28] , while for
$$
q \rightarrow 1 , h \rightarrow 1
$$
one recovers the content of the another example

in [28]. The case $h \rightarrow 0 $ is included

in the formulas of $q$-calculus of Thomae-Jackson easy to be
specified from [27] (see also up-date references there). For
Bernoulli- Taylor Formula (presented during $PTM$ - Convention
Lodz - $2002$) : contact [27] for its recent version.

\vspace{2mm}

Exercise: find the series expansion in of the $\partial
_{q,h}$-delta operator . Note and check that the polynomial
sequence $p_n(x) =\Big(x- [\frac{h}{1-q]}]\Big)^n , n\geq 0$ is
the Appell-Sheffer sequence for the $\partial _{q,h}$-delta
operator . This is consequently not recognized in [29] - see:
(5.6) - there.

\vspace{2mm}

\section{III.  Third   Contact Remark on $\tau$ -calculus and umbral
difference calculus }

In [29]  (page 2) the so called $\tau$ derivative was defined and
then the exposition of the idea of "$\tau$-calculus" follows with
an adaptation of standard  methods to obtain solutions  of the
eigenvalue and eigenvector  - this time -$\tau$-difference
equation of second order. (Compare $\partial_\tau$  of  [29] with
divided difference operator from  [30] )

\vspace{2mm}

Our first  $\tau$-experience and $\tau$-info  is: Inspection of
(5.7) formula and its neighborhood  from [29] leads to immediate
observation:

\vspace{2mm}

a]  if $X_n$ is standard monomial then  $\partial_\tau$ is a
$\partial_\psi$  derivative where $n_\psi =\frac {c_n}{c_{n-1}} ,
n>0 $, [3-5] .

b]  if $ X_n  = p_n(x)$  is a polynomial of n-th order then
$\left\{ {p_{n}}  \right\}_{n \ge 0} $, deg $p_{n} = n$ represents
$\psi$-Sheffer sequence and $ \partial_\tau = Q(\partial_{\psi})$
is the corresponding $\partial_\psi$-delta operator   so that
$Q(\partial_{\psi})p_n(x) = n_{\psi}p_{n-1}(x)$ [3-5]. This is not
recognized in [29] .

\vspace{2mm}

Our next $\tau$-experience and $\tau$-info is: the conclusive
illustrative {\bf Example 5.3} in [29] from 8 {\bf December} 2002
comes from the announced (p.18 and p.21 in [29]) A.Dobrogowska,
A.Odzijewicz \emph{Second $q$-difference equations solvable by
factorization method} .The calculations in ArXive preprint from
{\bf December} 2003 by A.Dobrogowska, A.Odzijewicz \emph{Second
$q$-difference equations solvable by factorization method}
$$ArXiv: math-ph/0312057 22 Dec 2003$$ rely  on calculations from
{\bf December} 2002 preprint [29] $$ArXiv: math-ph/020 8006.$$ On
the 14-th of January 2004 solemnly in public Professor Jan
Slawianowski had given the highest appraisal to this
investigation. He  compared its possible role to be played in the
history of physics with the Discovery of Max Planck from {\bf
December} 1900. Well. May be this is all because of {\bf
December}. Both preprints refer to difference Riccati equation.

\vspace{2mm}

One might be then perhaps curious about where $q$-Riccati equation
comes from according to the authors of [29].

\vspace{5mm}

{\bf Note:} before next section see: page 6 in  [29] [quotation]"
$q$-Riccati equation of [15]"[end of quotation]. Therefore note
also  that the $q$-Riccati equation has been considered in
productions of other authors earlier, for example [31-33]

{\bf IV. Fourth  Contact Remark on  $\psi$-calculus and
Viskov-Markowsky and Others -general  umbral difference calculus}

\vspace{2mm}

This section on perspectives is addressed  to those who feel Umbra
Power as well as to those who ignore  it  by  ....("umbra
dependent reasons" ? . ) What might be done  to extend the scope
of realized application of {\bf EFOC} = {\bf E}xtended {\bf
F}inite \textbf{O}perator {\bf C}alculus?

\vspace{2mm}

Certainly:\\
$\psi$-orthogonal OPS : orthogonality in Roman-Rota functional
formulation of umbral calculus [34] might be translated up to wish
into EFOC - language : see: [5,1]  for the beginings of
$\psi$-integration already in EFOC language - to be extended due
to inspiration by  $q,\tau , \psi$ - integrations , taking care of
roots coming from Jackson, Carlitz, Hann and so many Others ...

\vspace{2mm}

Certainly:\\
Certainly $\psi$-difference equations of higher order ... We may
first of all get more experienced and activate intuitions with
inspirations up to personal choice. Special up to my choice
sources of might be further $\psi$-inspirations are {\bf[35-50]} .

[35] by  Alphonse Magnus

{\bf contains }: ... first and second order difference equations -
orthogonal polynomials satisfying difference equations (Riccati
equations, continued fractions  etc. Pade, Chebyshev etc.)

{\bf contains }  many up-date relevant $q$-references and many
pages of substantial information - historical - included.

[36] by Ivan A. Dynnikov , Segey V. Smirnov :all six references:
principal substantial for the subject

[37] by  Mourad E. H. Ismail $q$-papers - most out of around 200
references ;{\bf contains} : $q$ -orthogonal polynomials,
Askey-Wilson divided difference operators,  inverses to the
Askey-Wilson operators .

[38,39,40,41]  $q$-Riccati equations,special functions,

$q$-series and related topics,$q$-commuting variables.

[42] by Haret C. Rosu , see also for history of $q$-difference
analogues

[43] by  I. Area, E. Godoy, F. Marcellán \textbf{contains}:
$q$-Jacobi polynomials; $q$-Laguerre/Wall polynomials...

[44,45,46,47] $q$-extension of the Hermite polynomials,coherent
pairs and orthogonal polynomials of a discrete variables,further
important $q$-analogs...

[48] by  Natig M. Atakishiyev, Alejandro Frank, and Kurt Bernardo
Wolf \textbf{contains}: Heisenberg $q$-algebra generators as
first-order difference operators, eigenstates  of the
$q$-oscillator Hamiltonian in terms of the $q$-1-Hermite
polynomials.... the measure for these $q$-oscillator states..

[49] by Brian D.  Taylor : Difference Equations via the Classical
Umbral Calculus

[50]  on  Umbral calculus and Quantum Mechanics

.......{\bf to be continued}.......

We end this incomplete Reference Survey attempted from  finite
operator calculus point of view with the link :
$$http://www.ms.uky.edu/~jrge/Rota/rota.html$$

It is the Gian-Carlo Rota and his grown up Children and Friends
-Family Link.


\begin{thebibliography}{99}
\parskip 0pt

\bibitem{1}
A.K.Kwasniewski {\it On Simple Characterizations of Sheffer
$\psi$-polynomials and Related Propositions of the Calculus of
Sequences} ,Bulletin de la Soc. des Sciences et de Lettres de
Lodz,  {\bf 52}  Ser. Rech. Deform.  36 (2002):45-65, ArXiv:
math.CO/0312397 .

\bibitem{2}
A. Odzijewicz  Commun. Math. Phys. {\bf 192}, 183 (1986).

\bibitem{3}
A.K.Kwa\'sniewski {\it Towards $\psi$-extension of Finite Operator
Calculus of Rota} Rep. Math. Phys. {\bf 48} (3), 305-342 (2001).
$$ArXiv: math.CO/0402078  2004$$

\bibitem{4}
A.K. Kwa\`sniewski {\it On extended finite operator calculus of
Rota and quantum groups}  Integral Transforms and Special
Functions {\bf 2}(4), 333 (2001)

\bibitem{5}
A.K.Kwa\`sniewski  {\it Main theorems of extended finite operator
calculus} Integral Transforms and Special Functions {\bf 14}, 333
(2003).

\bibitem{6}
M. Ward: {\it A calculus of sequences} Amer. J. Math. {\bf 58},
255-266 (1936).

\bibitem{7}
G. Markowsky {\it Differential operators and Theory of Binomial
Enumeration} J.Math.Anal.Appl. {\bf 63} (1978):145- 155

\bibitem{8}
Andrew P. Guinand {\it The umbral method: A survey of elementary
Mnemonic and Manipulative uses} American Math. Monthly {\bf 86}
(1979): 187- 195

\bibitem{9}
Brian D. Taylor {\it Difference Equations via the Classical Umbral
Calculus}  in Mathematical Essays in Honor of Gian-Carlo Rota,
Birkhauser, Boston, 1998.


\bibitem{10}
G.-C. Rota and B. D. Taylor, {\it The Classical Umbral Calculus}
SIAM J. Math. Anal. {\bf 25}, No.2 (1994):694-711.

\bibitem{11}
Ira M. Gessel {\it Applications of the classical umbral calculus
Algebra Universalis} {\bf 49} (2003): 397-434


\bibitem{12}
 Adams C. R. {\it Linear q-difference equations} Bull.Amer.Math.Soc.
 {\bf 37} (1931): 361-400

\bibitem{13}
T. Ernst    http://www.math.uu.se/~thomas/Lics.pdf

\bibitem{14}
$http://en2.wikipedia.org/wiki/Umbral_calculusUmbral$

\bibitem{15}
$http://www.google.pl/search?q=Blissard+symbolic+calculus$

\bibitem{16}
Di Bucchianico  , D. Loeb April {\it A Selected Survey of Umbral
Calculus} www.combinatorics.org/Surveys/ds3.pdf (2000)

\bibitem{17}
D. Loeb , G. C. Rota   {\it Recent Contributions to the calculus
of Finite Differences: a Survey}  Lecture Notes in Pure and Appl.
Math. {\bf 132} (1991):239-276, ArXiv: math.co/ 9502210 V1 9 Feb
1995

\bibitem{18}
[18] S. Roman     http://www.romanpress.com/me.htm

\bibitem{19}
$http://mathworld.wolfram.com/UmbralCalculus.html$

\bibitem{20}

R. Stanton $ http://www.math.umn.edu/~stanton/publist.html$

\bibitem{21}
NAVIMA  - group:
$$http://webs.uvigo.es/t10/navima/publications.html$$

\bibitem{22}
 Hahn W. " {\it Uber orthogonal Polynomen die q-Differenzengleichungen
genügen} Math. Nachr. Berlin {\bf 2} ({\bf 1949}): 4-34.

\bibitem{23}
J. Thomae   {\it Beirträge zur Theorie der durch die Heinesche
Reihe Darstellbaren Funktionen} J. reine angew. Math. {bf 70}
({\bf 1869}): 258-281


\bibitem{24}
F. H. Jackson  {\it q-Difference Equations} Amer.J.Math. {\bf 32}
({\bf 1910}): 305-314

\bibitem{25}
F. H. Jackson  {\it On q-definite integrals} Quart.J.Pure and
Appl. Math. {\bf 41} (1910): 193-203

\bibitem{26}
A. Odzijewicz et all. {\it Integrable multi-boson systems and
orthogonal polynomials} J.Phys.A : Math.Gen. {\bf 34} (2001):
4353-4376

\bibitem{27}
A.K.Kwa\`sniewski {\it Bernoulli-Taylor formula of psi-umbral
difference calculus}, preprint 2002/22 Faculty of Mathematics ,
Univesity of Lodz (2002) ArXiv: math.GM/0312401  (2003)

\bibitem{28}
O.V. Viskov {\it Trudy Matiematiczieskovo Instituta AN SSSR}
 {\bf 177} (1986):21 - 32

\bibitem{29}
A. Odzijewicz, T. Goli\`nski {\it Second order functional
equations} ArXiv: math-ph/020 8006 v2  8 Dec 2002


\bibitem{30}
Magnus A.  {\it New difference calculus and orthogonal
polynomials} ( 1997)
$$www.math.ucl.ac.be/~magnus/num3/m3011973.ps$$


\bibitem{31}
Gruenbaum F.A., Haine L. {\it On a $q$-analogue of Gauss equation
and some $q$-Riccati equations} American Mathematical Society,
Fields Inst. Commun.{\bf 14} 77-81 (1997).

\bibitem{32}
Jorge P. Zubelli {\it The Bispectral Problem, Rational Solutions
of the Master}  Amer. Math. Soc.Fields Inst. Commun.{\bf 14}(1997)
$$www.impa.br/~zubelli/PS/maiart.ps$$

\bibitem{33}
Michio Jimbo and Hidetaka Sakai {\it A $q$-Analogue of the six
Painlev`e equation} (1995)
$$www.kusm.kyoto-u.ac.jp/preprint/95/16.ps$$

\bibitem{34}
S. M. Roman     {\it The umbral calculus}   Academic Press  1984

\bibitem{35}
Alphonse Magnus, {\it Special topics in approximation theory}

New difference calculus and orthogonal polynomials. (1997)
$$
www.math.ucl.ac.be/~magnus/num3/m3011973.ps
$$
\bibitem{36}
Ivan A. Dynnikov , Segey V. Smirnov {\it Exactly solvable periodic
Darboux q-chains Russian Mathematical Surveys}, {\bf 57}
(2002):1218-1219. arXiv:math-ph/0207022 v1 18

\bibitem{37}
Mourad E. H. Ismail   www.math.usf.edu/~ismail/publications.html

\bibitem{38}
F. A. Grunbaum and L. Haine {\it On a q-analogue of Gauss equation
and some q-Riccati equations}, Special Functions, $q$-Series and
Related Topics", Edited by: Mourad E. H. Ismail , David R. Masson
and Mizan Rahman ,Amer.Math. Soc. Series:{\bf 14} (1997) : 77-81

\bibitem{39}
George Gasper {\it Elementary derivations of summation and
transformation formulas for $q$-series} (55--70);  ibid

\bibitem{40}
Erik Koelink [H. Tjerk Koelink] {\it Addition formulas for
q-special functions} (109--129); ibid


\bibitem{41}
Tom H. Koornwinder {\it Special functions and  $q$ -commuting
variables} (131--166);  ibid

\bibitem{42}
Haret C. Rosu {\it Short survey of Darboux transformations}
arXiv:quant-ph/9809056 v3 5 Oct 1999

\bibitem{43}
I. Area, E. Godoy, F. Marcellán, {\it $q$-Coherent pairs and
q-orthogonal polynomials}  Applied. Math. Comput.{\bf 128} (2002):
191-216.


\bibitem{44}
R. Álvarez-Nodarse, M. Atakishiyeva, N. M. Atakishiyev, {\it On a
$q$-extension of the Hermite polynomials $H_n(x)$ with the
continuous orthogonality property on R}, Bol. Soc. Mat. Mexicana
{\bf 8} (2002), 221-232.

\bibitem{45}
I. Area, E. Godoy, F. Marcellán, {\it D-coherent pairs and
orthogonal polynomials of a discrete variable}, Integral
Transforms Spec. Funct. {\bf 14} (2003): 31-57.


\bibitem{46}
Michio Jimbo and Hidetaka Sakai {\it A q-analog of the sixth
Painlev'e equation} (1995)
www.kusm.kyoto-u.ac.jp/preprint/95/16.ps


\bibitem{47}
R. Askey and SK Suslov  {\it The q-Harmonic Oscillator and an
Analog of the Charlier polynomials} arXiv:math.CA/9307206 v1 9 Jul
1993

\bibitem{48}
Natig M. Atakishiyev, Alejandro Frank, and Kurt Bernardo Wolf {\it
A simple difference realization of the Heisenberg q-algebra } J.
Math. Phys. {\bf 35}(7) 1994 : 3253-3260.

\bibitem{49}
Brian D. Taylor  {\it Difference Equations via the Classical
Umbral Calculus}  in Mathematical Essays in Honor of Gian-Carlo
Rota, Birkhauser, Boston, 1998.


\bibitem{50}
A.Dimakis, F. Muller-Hossei, T. Striker {\it Umbral calculus,
discretization, and Quantum Mechanics on a lattice} J. Phys. A:
Math. Gen. {\bf 29} (1996): 6861-6876 .

\end{thebibliography}
\end{document}